\documentclass[reqno]{amsart}
\usepackage{amsbsy,amssymb,amscd,amsfonts,latexsym,amstext,delarray,
 amsmath,graphicx, color}
 \usepackage[dvips]{epsfig}
\usepackage{hyperref}
\input xypic
\usepackage{srcltx}
\usepackage{mathrsfs}
\usepackage{float} 

\newtheorem{thm}{Theorem}[section]
\newtheorem{prop}[thm]{Proposition}
\newtheorem{cor}[thm]{Corollary}
\newtheorem{lem}[thm]{Lemma}

\def\Qer{{\rm Ker}}

\def\A{{\mathbb A}}

\def\C{{\mathbb C}}
\def\F{{\mathbb F}}

\def\P{{\mathbb P}}
\def\Q{{\mathbb Q}}
\def\R{{\mathbb R}}

\def\cE{{\mathcal E}}

\def\cL{{\mathcal L}}

\def\cS{{\mathcal S}}

\def\qqq{\,,\,~\forall}

\def\Lbd{\lambda}

\newcommand{\ie}{{\it i.e.\/}\ }
\newcommand{\eg}{{\it e.g.\/}\ }
\newcommand{\cf}{{\it cf.}}
\newcommand{\opcit}{{\it op.cit.\/}\ }

\def\sin{{{\rm sin}}}
\def\cos{{{\rm cos}}}

\def\dim{{\mbox{dim\,}}}
 \def\scal2{{\mathscr S}}
 
\def\fourier{\F}
\def\sr0{{\cS^{\rm ev}_0}}

\def\sar0{{\cS_0(\A_\Q)}}

\def\dom{{\rm Dom}}
\def\wmin{{W_{\rm min}}}
\def\wmax{{W_{\rm max}}}
\def\wsa{{W_{\rm sa}}}
\def\fourier{\F}
\def\fourierer{\fourier_{e_\R}}
\def\Fa{\F_{e_\R}}

\def\Dirac{{D \hspace{-7pt} \slash \,}}

\newcommand{\nil}[1]{}

\title[Prolate spheroidal operator]{Prolate spheroidal operator and Zeta}

\begin{document}
\author[Connes]{Alain Connes}
\author[Moscovici]{Henri Moscovici}
\address{A.~Connes: Coll\`ege de France \\
3, rue d'Ulm \\ Paris, F-75005 France\\
I.H.E.S. and The Ohio State University} \email{alain\@@connes.org}

\address{H.~Moscovici:
Department of mathematics, The Ohio State University, Columbus, OH 43210, USA}
\email{henri@math.ohio-state.edu}
\thanks{The work of H.~M. was partially
 supported by NSF award 1600541}

\begin{abstract}

In this paper we describe a remarkable new property of the self-adjoint extension $\wsa$ of the  prolate spheroidal operator introduced in \cite{college98},\cite{CMbook}.
The restriction of this operator to the interval $J$ whose
characteristic function commutes with it is well known, has discrete positive spectrum and is 
well understood \cite{Slepian0, Sl, Slepian, kat}.
What we have discovered is that the restriction of $\wsa$ to the complement of $J$ admits (besides a replica of the above positive spectrum) negative eigenvalues whose ultraviolet behavior reproduce that of the squares of zeros of the Riemann zeta function. 
Furthermore, their corresponding eigenfunctions belong to the Sonin space. This  feature
fits with the proof \cite{weilpos} of Weil's positivity at the archimedean place, which uses the compression of the scaling action to the Sonin space.  As a byproduct we construct an isospectral family of Dirac operators whose spectra have the same ultraviolet behavior as the zeros of the Riemann zeta function. 

\end{abstract}

\maketitle

\section{Introduction}

The prolate spheroidal wave functions play a key role in \cite{Co-zeta, weilpos,ccspectral} in relation with the Riemann zeta function. In all these applications they appear as eigenfunctions of the angle operator between two orthogonal projections in the Hilbert space $L^2(\R)^{\text{ev}}$ of even square integrable function on $\R$. These projections depend on a parameter $\lambda>0$, the projection $P_\lambda$ is given by the multiplication with the characteristic function of the interval  $[-\lambda,\lambda]\subset \R$. The projection $\widehat{P_\lambda}$ is its conjugate by the Fourier transform $\Fa$ which is 
 the unitary operator  in $L^2(\R)^{\text{ev}}$ defined by 
 $$
 \Fa(\xi)(y)=\int \xi(x)\exp(-2 \pi i x y) dx.
 $$
 In all the above applications of prolate spheroidal wave functions the miraculous existence, discovered by the Bell Labs group \cite{Slepian0, Sl, Slepian}, of a differential operator $W_\lambda$ commuting with the angle operator, plays only an auxiliary role. 
 In the present paper we uncover another ``miracle": a careful study of the natural self-adjoint extension of $W_\lambda$ introduced in \cite[Lemma 6]{college98} 
(see also \cite[\S 3.3]{CMbook}) to $L^2(\R)$ shows that it still has discrete spectrum 
and that its negative eigenvalues  reproduce the ultraviolet behavior  of the squares of zeros of the Riemann zeta function. In a similar way the positive spectrum corresponds, in the ultraviolet regime,  to the trivial zeros. This coincidence holds for two values  $\lambda=1$ and $\lambda=\sqrt 2$. The conceptual reason for this coincidence is the link between the operator  
\begin{equation}\label{prol1}
	(W_\lambda\xi)(x)=-\partial_x (\lambda^2 -x^2) \partial_x \xi(x)+ (2\pi \lambda)^2 x^2 \xi(x)
\end{equation}
and the square of the scaling operator $S:=x\partial_x$. In \cite{weilpos} the compression of $f(S)$ to Sonin's space (for $\lambda=1$) was shown to be the root of Weil's positivity at the archimedean place on test functions with support in the interval $[2^{-1/2},2^{1/2}]$, but since Sonin's space is not preserved by scaling, one could not restrict scaling to this space. It turns out that $W_\lambda$ commutes with the orthogonal projection on Sonin's space.
Thus one can restrict $W_\lambda$ to Sonin's space and the ultraviolet spectral similarity with the squares of non-trivial zeros of zeta suggests that one has spectrally captured the contribution of the archimedean place to the mysterious zeta spectrum. In fact using the Darboux process we construct a Dirac square root of $W_\lambda$ depending on a deformation parameter, and whose spectrum has the same ultraviolet behavior as the zeros of the Riemann zeta function. \newline

Our paper is organized as follows: In Section \ref{sectsa} we show that there exists a unique selfadjoint extension $\wsa$ of the symmetric operator $\wmin$ defined on 
Schwartz space $\cS(\R)$ by \eqref{prol1}. Moreover $\wsa$ commutes with Fourier transform 
and has discrete spectrum unbounded in both directions. In Section \ref{sectsonin} we show that the eigenvectors for negative eigenvalues of $\wsa$ belong to Sonin's space. In Section \ref{sectsemiclassical} we compute the semiclassical approximation to the number of negative eigenvalues of $\wsa$ whose absolute value is less than $E^2$. 
In Section \ref{sectdirac}, we use the Darboux method combined with solutions of a Riccati 
equation to construct an isospectral family of Dirac operators $\Dirac$ whose squares are direct sums of two copies of $\wsa$. In Section \ref{sectroot2} we specialize to the case $\lambda=\sqrt 2$   and show that the operator $2 \Dirac$ has discrete simple spectrum contained in $\R\cup i\R$ with imaginary eigenvalues  symmetric under complex conjugation and  counting function $N(E)$ (counting those of positive imaginary part less than $E$) fulfilling the same as the  Riemann formula 
\begin{equation}\label{neintro}
N(E)\sim \frac{E}{2 \pi }\left(\log \left(\frac{E}{2 \pi }\right)-1\right)+O(1)
\end{equation}
We also show the numerical evidence for the ultraviolet spectral similarity between the eigenvalues of $2\Dirac$ and the zeros of the Riemann zeta function. Lastly, Section \ref{sectrems} contains more speculative final remarks, in particular on a natural 
two-dimensional black hole geometry intrinsically related to the operator $2\Dirac$.

\section{The selfadjoint prolate wave operator}\label{sectsa}

The prolate spheroidal operator \eqref{prol1} is an operator of Sturm-Liouville type, 
\begin{align}\label{prol2}
\begin{split}
	(W_\lambda\xi)(x) &= -\partial_x \big(p(x) \partial_x \xi(x)\big)+ q(x)\xi(x) , \qquad x \in \R \\
\text{where} &\qquad  p(x) = \Lbd^2-x^2, \quad  q(x) = (2\pi \lambda)^2 x^2 ,
\end{split}
\end{align}
but having two interior singular points it is not directly treatable by the usual 
Sturm-Liouville theory. However
its restrictions to each of the intervals $(-\infty, -\Lbd)$, $( -\Lbd, \Lbd)$ and 
$(\Lbd, \infty)$ are standard, in fact quasi-regular, Sturm-Liouville operators.

\

Henceforth $W_\Lbd$  will be simply denoted $W$ whenever $\Lbd$ is a general parameter. 
To begin with, we regard $W$ as an unbounded operator 
on $L^2(\R)$ with core the Schwartz space $\cS(\R)$.
As such, $W$ is real, symmetric and invariant under the parity exchange $x\mapsto -x$.  
These features are inherited by its closure in the graph norm  $\wmin$, as well as by  
 $\wmax=\wmin^*$, the latter having domain
\begin{equation}\label{op1adjoint}
	\dom(\wmax)=\{\xi \in L^2(\R)\mid W \xi \in L^2(\R) \} ,
	\end{equation}
with $W \xi$ viewed as a tempered distribution. 
In addition $W$ has the remarkable property of commuting with the Fourier transform
 \begin{equation}\label{fourier}
	\fourierer(f)(y):=\int _{-\infty}^\infty f(x) \exp{(-2 \pi i x y)}dx .
\end{equation}
Since both the Schwartz space $\cS(\R)$ and its dual are
globally invariant under the Fourier transform, the domains $\dom \wmin$ and $\dom \wmax$ 
are invariant too, therefore both $\wmin$ and $\wmax$ commute with $\fourierer$.
  
\begin{lem}\label{Wdefic} The deficiency indices of $\wmin$ are  $(4,4)$.	
\end{lem}

\proof Any $\xi \in \dom(\wmax)$ satisfying $W\xi=\pm i\xi$ is a piecewise real analytic 
function and is uniquely specified by six parameters in the complement of the two regular 
singular points $\pm \lambda$. The known form of the solutions (\cf. \cite{ramisrecent})
together with the fact that $W \xi \in L^2(\R)$ imply that the logarithmic singularities 
of $\xi$ on the left and
the right of $\pm\lambda$ have to match. This reduces the number of parameters to $4$.
Conversely, since all $4$ singular points are LC (limit circle case), any solution of $W\xi=\pm i\xi$ belongs
to $\dom(\wmax)$, hence $\dim \Qer(\wmax \pm i I) =4$.
\endproof

\begin{lem}\label{regsinglem} Let $\xi \in \dom \wmax$ and denote $a=\pm \lambda$.
 The distribution $p(x)\partial_x\xi$ coincides with a continuous function $f$ in a neighborhood of $a$ and the evaluation map $L(\xi):=f(a)$ 
defines a non-zero continuous linear form on $\dom \wmax$ which  vanishes on the closed subspace $\dom\wmin$.	\newline
\end{lem}

\proof  Let $V=[b,c]$ be a compact interval neighborhood of $a=\pm \lambda$ where $a$ is the only zero of $p(x)$. Let $\psi$ be the distribution $\psi=p(x)\partial_x\xi(x)$, one has by definition,
$$
\langle \psi \mid \phi\rangle= -\int_\R \xi(x)\partial_x (p(x)\phi(x))dx \qqq \phi \in \cS(\R)
$$
Let $\eta=\wmax \xi$, one has by definition,
$$
\langle \eta \mid \phi\rangle=\langle \xi \mid \wmin\phi\rangle =\int_\R \xi(x)\left(\partial_x (p(x)\partial_x\phi(x))+q(x)\phi(x)\right)dx \qqq \phi \in \cS(\R)
$$ 
Let $\xi_1\in L^2(\R)$ coincide with $q(x)\xi(x)$ on $V$. Then for any smooth function
$\phi$ with support in $V$,
$$
\langle \psi \mid \partial_x\phi\rangle= -\int_\R \xi(x)\partial_x (p(x)\partial_x\phi(x))dx=
\langle \xi_1-\eta \mid \phi\rangle
$$  
The restriction of $\xi_1-\eta$ to $V$ belongs to $L^2(V)\subset L^1(V)$ and the function $f_1(x)=-\int_b^x(\xi_1-\eta)(t)dt$ is continuous and fulfills 
$$
\int_V f_1(x)\partial_x\phi(x)dx=\langle \xi_1-\eta \mid \phi\rangle
$$
It follows that $\langle \psi-f_1 \mid \partial_x\phi\rangle=0$ for all smooth functions $\phi$ with support in $V$ and choosing a positive smooth function $\phi_1$ with support in $V$ and integral $1$, one obtains 
$$
\langle \psi \mid \phi\rangle=\langle f_1+ s\mid \phi\rangle \qqq \phi \in C_c^{\infty}(V), \ s=\langle (\psi-f_1)\mid \phi_1\rangle .
$$
Thus the distribution $p(x)\partial_x\xi$ coincides with the function $f(x):=f_1(x)+s$ on $V$. One has 
$$
f(a)=s+ f_1(a)=\langle (\psi-f_1) \mid \phi_1\rangle-\int_b^a(\xi_1-\eta)(x)dx
$$ 
Moreover $\langle \psi \mid \phi_1\rangle=\int p(x)\partial_x\xi(x)\phi_1(x)dx=-\int \xi(x) \partial_x(p(x)\phi_1(x))dx=\langle \xi \mid \eta_1\rangle$ where $\eta_1\in C_c^{\infty}(V)$. One has also, 
$$
-\langle f_1 \mid \phi_1\rangle-\int_b^a(\xi_1-\eta)(x)dx=\int_b^c\int_b^x(\xi_1-\eta)(t)\phi_1(x)dtdx-\int_b^a(\xi_1-\eta)(x)dx$$ $$=
\langle \xi \mid \eta_2\rangle+\langle \eta \mid \eta_3\rangle
$$
where the vectors $\eta_j\in L^2(\R)$. Thus the linear form $L(\xi):=f(a)$ is continuous in the graph norm of $\dom\wmax$. For $\xi\in \cS(\R)$ the distribution $\psi=p(x)\partial_x\xi(x)$ is a function vanishing at $x=a$ and thus $L(\xi)=0$. By the density of $\cS(\R)$ in $\dom \wmin$ for the graph norm, it follows that $L$ vanishes on the closed subspace $\dom\wmin$.	
\endproof  

Let $P_\lambda$ be the cutoff projection associated to the interval $[-\lambda,\lambda]$, \ie the multiplication operator by the characteristic function $1_{[-\lambda,\lambda]}$, and let 
$\widehat P_\lambda = \fourierer P_\lambda \fourierer^{-1}$ denote its conjugate by the Fourier transform.

\begin{lem}\label{WP} 
If $\xi\in \dom \wmin$  then $P_\lambda \xi \in \dom \wmax$ and 
$W P_\lambda \xi= P_\lambda W \xi$. The same holds with respect to $\widehat P_\lambda$.
\end{lem}

\proof  
Let $f\in C^{\infty}(V)$ where $V$ is a neighborhood of the interval $[-\lambda,\lambda]$.
 Then $P_\lambda f\in \dom \wmax$ and viewing $W(P_\lambda f)$ as distribution one gets,
 for any $ \phi \in \cS(\R)$
\begin{align*}
\langle W(P_\lambda f), \phi\rangle &=\int_{-\lambda}^\lambda f(x)(W\phi)(x)dx=\int_{-\lambda}^\lambda -f(x)\partial_x (\lambda^2 -x^2) \partial_x \phi(x)dx \\
&+\int_{-\lambda}^\lambda f(x) (2\pi \lambda)^2 x^2 \phi(x)dx .
\end{align*}
Using  twice integration by parts, together with
 the fact that $(\lambda^2 -x^2)\phi'(x)$ and $(\lambda^2 -x^2)f'(x)$ vanish on the boundary, 
 one obtains
\begin{align*}
\langle W(P_\lambda f), \phi\rangle&=\int_{-\lambda}^\lambda f'(x)((\lambda^2 -x^2)\phi')(x)dx+\int_{-\lambda}^\lambda f(x) (2\pi \lambda)^2 x^2 \phi(x)dx \\
 &=-\int_{-\lambda}^\lambda (\partial_x((\lambda^2 -x^2)f'(x)))\phi(x)dx+\int_{-\lambda}^\lambda f(x) (2\pi \lambda)^2 x^2 \phi(x)dx\\
&=\int_{-\lambda}^\lambda Wf(x)\phi(x)dx ,
\end{align*}
 which shows that $W(P_\lambda f)=P_\lambda Wf$.  
  In particular the same is true for any $f\in \cS(\R)$, and by the density 
 of $\cS(\R)$ in $\dom \wmin$ for the graph norm it follows that 
\begin{align*}
\xi\in \dom \wmin \, \Longrightarrow \, P_\lambda \xi \in \dom \wmax \quad \text{and} \quad 
\wmax P_\lambda \xi= P_\lambda W \xi .
\end{align*}
The claim now follows from the fact that $W$ commutes with $\fourierer$.
  \endproof

The selfadjoint extensions of $\wmin$ are
parametrized by self-orthogonal subspaces of $\cE:= \dom(\wmax)/\dom(\wmin)$
with respect to the anti-symmetric sesquilinear form given by the pairing
\begin{equation} \label{O-def}
 \Omega(\xi,\eta) := \,\frac 1 i \Bigl( \langle \wmax \xi \mid \eta \rangle -  
 \langle \xi \mid \wmax \eta \rangle \Bigr), 
 \qquad \xi, \eta \in  \dom(\wmax)
  \end{equation}
which descends to a  non-degenerate form on $\cE$. 

The $\Omega$-pairing can be expressed in terms of boundary values as usual.
One starts with the Lagrange identity
 \begin{equation}
 \frac{d}{dx} [\xi, \eta] \, = \, \xi \, W\eta - \eta \, W \xi ,  
 \end{equation}
 where $\xi, \eta \in C^1(\R) \cap \dom\wmax$, and
  \begin{equation}
[\xi, \eta] := p \left(\xi\frac{d\eta}{dx} -  \eta\frac{d\xi}{dx}\right)  , \qquad p(x)= \Lbd^2 - x^2 ,
 \end{equation}
 is the (generalized) Wronskian.
 By integrating it on compact subintervals $[a, b] \subset \R \setminus \{\pm \lambda\}$ 
one obtains Green's formula
 \begin{align} \label{Lag}
 \int_a^b \left( W(\xi) \bar{\eta} - \xi W(\bar{\eta})\right)(x) dx 
 = [\xi, \bar{\eta}]\vert_a^b \, := \lim_{x \to b} [\xi, \bar{\eta}](x) - \lim_{x \to a} [\xi, \bar{\eta}](x) .
  \end{align}
Passage to the lateral limits towards the endpoints of the three subintervals partitioning
 $\R$ extends this identity to the whole real line, allowing to express $\Omega$
 in terms of Lagrange brackets as follows:  
\begin{align} \label{O-Lag}
i \Omega(\xi,\eta) =  [\xi, \bar{\eta}]\vert_{-\infty}^{-\Lbd} + [\xi, \bar{\eta}]\vert_{-\Lbd}^{\Lbd} +
 [\xi, \bar{\eta}]\vert_\Lbd^\infty 
\end{align} 
for all pairs $\xi, \eta \in \dom\wmax$.
 
\

Since $W$ is invariant under parity exchange, it preserves the orthogonal decomposition  
$L^2(\R)=L_+^2(\R) \oplus L_-^2(\R)$ into  even, resp. odd functions, which in turn induces
corresponding splittings $W=W^+ \oplus W^-$, $\Omega=\Omega_+. \oplus \Omega_-$ 
and $\cE= \cE_+ \oplus \cE_-$. Note also that $\cE_\pm$ are invariant under Fourier 
transform.

The following auxiliary lemma will be used in the ensuing discussion. 
\begin{lem}\label{atinfty} 
$(i)$~Let $f(x)=\frac 12\log ((\lambda^2-x^2)^{-2})$	viewed as a tempered distribution. Then the Fourier transform $\fourierer f$ is a distribution which coincides outside $0$ with the function 
$$\widetilde f (y)= \frac{\cos(2\pi \lambda y)}{ \vert y\vert} .$$
$(ii)$~Let $1_I$ be the characteristic function of the interval $I=[-\lambda,\lambda]$ then $$\fourierer 1_I(y)=\frac{\sin(2\pi \lambda y)}{\pi y} .$$ 
\end{lem}
\proof 
$(i)$~One has $f(x)=\frac 12\log ((\lambda-x)^{-2})+\frac 12\log((\lambda+x)^{-2})$, thus we start by computing  the Fourier transform of the distribution $\ell = -\log (x^2)$. 
One has $x\partial_x \ell=-2$. Thus one gets  $\partial_y y\widehat{\ell}=2$. Therefore
$y\widehat{\ell}$ is equal to ${\rm sign}(y)$ and $\widehat{\ell}$ is the Weil principal value $1/\vert y\vert$. Translation of the variable means multiplication by an imaginary exponential in Fourier and this gives the required equality. \newline
$(ii)$~One has $\partial_x 1_I=\delta_{-\lambda}-\delta_{\lambda}$ and in general $\fourierer f(y)=2\pi i y\, \fourierer f(y)$.\newline\endproof

We now proceed to construct a basis of $\cE$. 
First, for $\cE_+$ we pick an even function  $\alpha_+ \in C_c^\infty(\R)$ such that 
$\alpha_+(x)= \log |\lambda^2 -x^2|$
for $x \in [\frac 34 \Lbd, \frac 54 \Lbd]$ and with support in $(\frac 12 \Lbd, \frac 32 \Lbd)$. 
Then  we take $\beta_+(x) = 1_I$, the characteristic function of the interval 
$I=[-\lambda,\lambda]$, which belongs to $P_\lambda \cS(\R)$ and hence to $\dom\wmax$.
Next for $\cE_-$ we let $\alpha_-(x):=x \alpha_+(x)$ and $\beta_-(x):=x \beta_+(x)$.

\begin{lem}
The quadruplet $\{\alpha_\pm,  \beta_\pm, \widehat\alpha_\pm,
\widehat\beta_\pm \}$ forms a basis of  $\cE_\pm$. 
\end{lem}
\proof One checks using the expression \eqref{O-Lag} of the $\Omega$-pairing together
with Lemma \ref{atinfty} that the matrix representation of $\Omega_+$ with respect to the
given quadruplet has a single nonzero entry in each row and column.

In the odd case we note that  on the one hand 
 $ [\alpha_-, \beta_-](x) = x^2 [\alpha_+, \beta_+](x)$, and on the other hand
the derivatives involved in their Fourier transforms 
$\widehat\alpha_-(x) =  \frac{i}{2\pi}\partial_x\widehat\alpha_+(x)$, resp. 
$\widehat\beta_-(x) = \frac{i}{2\pi}\partial_x\widehat\beta_+(x)$.
exchange the two functions $cos$ and $sin$ in the leading terms at infinity. With this observation the calculation becomes similar
 to that for the even case, and so is the result.
\endproof

\

The $\Omega$-pairings with the above basis elements yield boundary conditions 
 of Sturm-Liouville type. Denoting, for $\xi \in \dom(W^\pm_{\max})$,
\begin{align}  \label{BdL}
\begin{split}
\L_{\alpha_\pm}(\xi) :=&\, i \Omega_\pm(\xi, \alpha_+) , \qquad
\L_{\widehat\alpha_\pm}(\xi) :=\, i\Omega_\pm(\xi, \widehat\alpha_+)  , \\
\L_{\beta_\pm}(\xi) :=&\, i \Omega_\pm(\xi, \beta_+)  ,  \qquad
\L_{\widehat\beta_\pm}(\xi) = \, i \Omega_\pm(\xi, \widehat\beta_-) ,
 \end{split}
\end{align} 
the minimal domains are characterized in these terms as being the intersection
\begin{align} \label{dmin}
\dom(W^\pm_{\min}) = \Qer\, \L_{\alpha_\pm} \cap \Qer \, \L_{\beta_\pm} \cap
\Qer\, \L_{\widehat\alpha_\pm} \cap \Qer \, \L_{\widehat\beta_\pm} 
 \end{align} 
and the induced functionals on $\cE_\pm=\dom(W^+_{\max})/\dom(W^\pm_{\min})$ 
form a basis of $\cE_\pm^*$. 

By straightforward calculation, using the fact that one can always
restrict the computation to $\R^+$, one
obtains explicit expressions for the boundary functionals. Up to a nonzero
constant factor they are as follows.  
In the even case,
\begin{align} \label{L-even}
\begin{split}
\L_{\alpha_+}(\xi) = &\, \lim_{x \nearrow \Lbd}\left((x-\Lbd)\log(\Lbd-x) \partial_x\xi(x) 
- \xi(x)\right) \\
 &\qquad \qquad - \lim_{x \searrow \Lbd}\left((x-\Lbd)\log(x-\Lbd) \partial_x\xi(x) 
 - \xi(x)\right) ; \\
\L_{\beta_+}(\xi) :=&\,   \lim_{x \nearrow \Lbd} \left((\Lbd -x)\partial_x\xi(x)\right)
 =  \lim_{x \searrow \Lbd}
\left((\Lbd -x)\partial_x\xi(x)\right) ; \\   
\L_{\widehat\alpha_+}(\xi) := & \,  \frac 2 \pi \lim_{x \to \infty}
\left(x \cos(2\pi\Lbd x)\partial_x\xi(x) +\bigl( 2\pi \Lbd x\sin(2\pi\Lbd x) + \cos(2\pi\Lbd x)\bigr)
\xi(x) \right) ;\\
\L_{\widehat\beta_+}(\xi) :=& -\frac 2 \pi \lim_{x \to \infty}
\left(x \sin(2\pi\Lbd x)\partial_x\xi(x)-\bigl( 2\pi \Lbd x\cos(2\pi\Lbd x) - \sin(2\pi\Lbd x)\bigr)
\xi(x) \right) .
\end{split}
\end{align} 
We note that the existence of the limit defining  $\L_{\beta_+}(\xi)$, \ie the equality of the
lateral limits, is ensured by Lemma \ref{regsinglem}.
 
Similar formulas define the functionals $\L_{\alpha_-}, \L_{\beta_-}, \L_{\widehat\alpha_-},
\L_{\widehat\beta_-}$ in the odd case.

 \

 Since both $\dom(\wmin)$ and $\dom(\wmax)$, as well as the symplectic form $\Omega$,
 are globally invariant under the Fourier transform, the quotient inherits
 induced transformations $f^\pm_{e_\R}: \cE_\pm \to \cE_\pm$ which relates the
boundary functionals as follows:
\begin{align}
 \L_{\widehat\beta_\pm}  = \L_{\beta_\pm}\circ f_{e_\R} \quad \text{and} \quad
  \L_{\widehat\alpha_\pm}  = \L_{\alpha_\pm} \circ f_{e_\R} .
\end{align}  
This association gives rise to two distinguished self-orthogonal subspaces, namely
\begin{align}
 \cL_{\beta}  = \bigcap_\pm \Qer \L_{\beta_\pm} \cap  \bigcap_\pm \Qer \L_{\widehat\beta_\pm} 
 \quad  \text{and} \quad 
 \cL_{\alpha}  = \bigcap_\pm \Qer \L_{\alpha_\pm} \cap  \bigcap_\pm 
 \Qer \L_{\widehat\alpha_\pm} 
\end{align}

\ 

{\bf Definition.} We denote by $\wsa$ the restriction of the operator $\wmax$ to the subspace 
$ \cL_{\beta}  = 
\bigcap_\pm \Qer \L_{\beta_\pm} \cap  \bigcap_\pm \Qer \L_{\widehat\beta_\pm} $.
Explicitly, its domain $\dom\wsa$ consists of the elements 
 $\xi \in \dom(W_{\max})$ satisfying the following boundary conditions:
 \begin{align} \label{I-1}
\lim_{x \to \pm\Lbd}(\Lbd^2 - x^2)\partial_x\xi(x) = 0  ,
\end{align} 
and at $\pm \infty$, writing $\xi=\xi^+ + \xi^-$ with $\xi^\pm \in \dom(W^\pm_{\max})$, 
 \begin{align} \label{l-ev}
& \lim_{x \to \pm\infty}  
\left(x \sin(2\pi\Lbd x)\partial_x\xi^+(x) -\bigl( 2\pi \Lbd x\cos(2\pi\Lbd x) - \sin(2\pi\Lbd x)\bigr)
\xi^+(x)\right) =0 , \\
& \lim_{x \to \pm\infty} \label{l-odd} 
\left(x \cos(2\pi\Lbd x)\partial_x\xi^-(x) +\bigl( 2\pi \Lbd x\sin(2\pi\Lbd x) + \cos(2\pi\Lbd x)\bigr)
\xi^-(x)\right) =0 .
\end{align} 
 
 \
 
We are now in a position to establish the main result of this section.

 \begin{thm}\label{thmwsa}  
  $(i)$~$\wsa$ is selfadjoint and commutes with the Fourier transform. \newline 
  $(ii)$~$\wsa$ commutes with the projections $P_\lambda$ and $\widehat P_\lambda$. \newline
 $(iii)$~$\wsa$ is the only selfadjoint extension of $\wmin$ commuting with 
 $P_\lambda$ and $\widehat P_\lambda$.
  $(iv)$~The spectrum of $\wsa$ is discrete and unbounded on both sides.
 \end{thm}

 \proof $(i)$~$\wsa$ is selfadjoint by construction, and its
domain $\cL_{\beta}$ is invariant under the Fourier transform also 
by construction.  
 
 $(ii)$~Since $\dom \wmin $ is given by \eqref{dmin}.
every element of $ \cL_{\beta}$ is a linear combination of an element $\xi \in \dom \wmin$ and
the $4$ vectors $\beta_\pm, \widehat \beta_\pm$ of Lemma \ref{Wdefic}. Each $\beta_\pm$ 
is of the form $P_\lambda f_\pm $ with $f_\pm$ smooth with compact support and thus one has, using Lemma \ref{WP},
$$
P_\lambda \beta_\pm=\beta_\pm\in \cS , \ \ \wsa P_\lambda \beta_\pm=\wsa P_\lambda f_\pm= P_\lambda Wf_\pm ,
$$
which shows that $\wsa P_\lambda \beta_\pm=P_\lambda \wsa P_\lambda \beta_\pm=P_\lambda \wsa \beta_\pm$ giving the required commutation for the $\beta_\pm$.

$(iii)$~The domain of a selfadjoint extension of $\wmin$ commuting with $P_\lambda$ and $\widehat P_\lambda$ must be contained in $\dom \wmax$ and also
contain both $P_\lambda\cS(\R)$ and $\widehat P_\lambda\cS(\R)$. Thus it must contain 
$\cL_{\beta}$, and cannot be larger due to self-adjointness.

$(iv)$ The operators $P_\lambda \wsa$ and $(I-P_\lambda) \wsa$ are selfadjoint
on $(-\lambda, \lambda)$, resp. on $(-\infty, \lambda) \sqcup (-\lambda, \infty)$, 
and thus covered by standard results in Sturm-Liouville theory (\cf \cite{naimark}, \cite{weid}, \cite{kat}). 
Indeed all four endpoints are \textit{limit circle case} in the Weyl classification 
(see \eg \cite[\S\S 5-6]{weid} for relevant definitions and properties), which
can be easily checked by using explicit bases of formal solutions 
for $W\xi -\mu \xi =0$, $\mu\in \C$,  around each singular point 
(\cf  \eg \cite[\S 2]{ramisrecent}).
The endpoints $\pm \lambda$ are LCNO (non-oscillatory limit circle), while
$\pm \infty$ are LCO (oscillatory limit circle) since the prolate spheroidal wave functions 
(which provide principal solutions around  $\pm \lambda$) have
infinitely many zeros in the neighborhood of  $\pm \infty$ (\cf \cite{Slepian0, Sl, Slepian}).
By well-known results (\cf \eg \cite{naimark} page 90, \cite{zettl}) it follows that
both $P_\lambda \wsa$ and $(I-P_\lambda) \wsa$ have discrete spectrum
and that the spectrum of the latter is unbounded on both sides. 
\endproof

\

 \begin{cor} \label{eigenfct} If $\phi$ is an eigenfunction of $\wsa^\pm$ then
\begin{itemize}
\item[{\bf (i)}] $\phi$ is regular on $[\Lbd, \Lbd+\epsilon)$ and on  $(\Lbd -\epsilon, \Lbd]$
for some $\epsilon > 0$, with a possible discontinuity at $\Lbd$;

\item[{\bf (ii)}] the leading term of the asymptotic expansion of $\phi$ at $\infty$ is 
proportional to $\frac{\sin(2\pi\Lbd x)}{x}$ if $\phi$ is even and to 
$\frac{\cos(2\pi\Lbd x)}{x}$  if $\phi$ is odd.
\end{itemize} 
 \end{cor}

\proof This follows from the above characterization \eqref{I-1}, \eqref{l-ev}, \eqref{l-odd}
of the domain of $\wsa$ combined with
the known bases of formal solutions for the equation $ W\xi = \mu \xi$, $\mu \in \R$, 
around $\pm \lambda$ and $\pm \infty$ (\cf \cite{ramisrecent}). 
\endproof


\section{Sonin space and negative eigenvalues}\label{sectsonin}

We translate the requirement that the Fourier transform $\fourierer f$ of an $f\in \dom \wmax$ has no logarithmic singularity at the singular points into a condition on the asymptotic behavior of $f$ at $\infty$. For simplicity we only deal with even functions, and for notational convenience
take $\Lbd=1$.

We can then find the asymptotic expansion at $\infty$ using the boundary condition that the leading term there is $\frac{\sin(2\pi \lambda y)}{ y}$. We take for simplicity $\lambda =1$ and use \cite{ramisrecent} to get for the tentative eigenvector for eigenvalue $\mu$ the expansion at $\infty$
$$
\xi_\mu(x)\sim \frac{\sin(2\pi  x)}{ x}+\frac{(\mu-4\pi^2)\cos (2 \pi  x)}{4\pi x^2}+\frac{-\mu ^2+8 \pi ^2 \mu +2 \mu -16 \pi ^4+8 \pi ^2}{32 \pi ^2 x^3}\sin(2\pi  x) +0(x^{-4})
$$
In fact as shown in Proposition 14 of \cite{ramisrecent}, the coefficients of this expansion are directly related to the coefficients of the expansion of the finite solution at $\lambda$ and taking for simplicity $\lambda =1$, if the latter is of the form 
$$
f_\mu(x)=\sum U_n(\mu) (x-1)^n, \ \ U_0(\mu)=1, \ U_1(\mu)=\frac{\mu-4\pi^2 }{2}
$$
$$
U_2(\mu)=\frac{\mu ^2-8 \pi ^2 \mu -2 \mu +16 \pi ^4-8 \pi ^2}{16}, \ \ldots 
$$
then the asymptotic series at infinity which governs the solution which has leading term in $\exp(-2 \pi i x)/x$ is equal to $v(x)\exp(-2 \pi i x)/x$ where
$$
v(x)\sim \sum n! U_n(\mu)(2\pi i x)^{-n}
$$
When one applies the Borel summation to this series the first step is to replace it by its Borel transform which is, up to normalization, 
$$
B(y):=\sum U_n(\mu) y^n
$$
and is related to $v(x)$ by $\int_0^{\infty } t^n \exp (-zt) \, dt=z^{-n-1} \Gamma (n+1)$ \ie the Laplace transform
$$
\frac{v(x)}{2\pi i x}=\int_0^{\infty} \exp(-2\pi i xt)B(t)dt
$$
\begin{lem} For any $\mu \in \R$ the asymptotic expansion of the unique solution $\xi_\mu$ which at $\infty$ is asymptotically  $\sim -\frac{\sin(2\pi  x)}{ \pi x}$ is Borel summable and is equal to the Fourier transform of the unique even solution $\phi_\mu$ which is zero on $[-1,1]$ and agrees with $f_\mu(x)$ for $x>1$.
\end{lem}
\proof One has the equality 
$$
\frac{v(x)}{2\pi i x}=\int_0^{\infty} \exp(-2\pi i xt)B(t)dt
=\int_0^{\infty} \exp(-2\pi i xt)f_\mu(t+1)dt=
$$
$$
=\int_1^{\infty} \exp(-2\pi i x(y-1))f_\mu(y)dy
$$
Thus one gets
$$
v(x)\exp(-2 \pi i x)/(2\pi i x)= \int_1^{\infty} \exp(-2\pi i xy)f_\mu(y)dy
$$
The function $\phi_\mu$ is even and vanishes on $[-1,1]$ so
$$
\int_{-\infty}^\infty \exp(-2\pi i xy)\phi_\mu(y)dy=
 \int_1^{\infty} \exp(-2\pi i xy)f_\mu(y)dy
+\overline{\int_1^{\infty} \exp(-2\pi i xy)f_\mu(y)dy}=
$$
$$
=v(x)\exp(-2 \pi i x)/(2\pi i x)+\overline{v(x)\exp(-2 \pi i x)/(2\pi i x)}
$$
Now these two terms are asymptotic solutions since $\mu$ is real and $v(x)\exp(-2 \pi i x)/(2\pi i x)$ is an asymptotic solution. Moreover the leading behavior at $\infty$ is in 
$$
\exp(-2 \pi i x)/(2\pi i x)-\exp(2 \pi i x)/(2\pi i x)=-\frac{\sin(2 \pi x)}{\pi x}
$$
Thus it follows that the Fourier transform  $\fourierer \phi_\mu=\xi_\mu$. 
\endproof 

 \
 
 \begin{cor} \label{soninfct}
 With the above notation, assume $\mu$ is a negative eigenvalue. Then $\phi_\mu$ 
 belongs to the Sonin space.
\end{cor}
\proof In fact Sonin's space is the  orthogonal  of the eigenspaces of $\wsa$ associated to the classical prolate functions and their Fourier transforms. \endproof
We should note that at this point we do not claim (although this is supported by numerical evidence)  that all eigenvalues of the restriction of $\wsa$ to  Sonin's space are negative, however there could be only finitely many exceptions. 
 
\section{Semiclassical approximation and counting function} \label{sectsemiclassical}

In this section we use the semiclassical estimate for the function counting the number of eigenvalues 
and investigate the negative eigenvalues of the operator $\wsa$. We consider the classical Hamiltonian
\begin{equation}\label{hamilt}
H_\lambda(p,q)=(p^2-\lambda^2)(q^2-\lambda^2)
\end{equation}
and use it as a semiclassical approximation of $\wsa$ via the formal relation 
\begin{equation}\label{hamilt1}
W_\lambda \sim -4\pi^2 H_\lambda+4\pi^2 \lambda^4
\end{equation}
using the correspondence $q\to x$ and $p\to \frac{1}{2\pi i}\partial_x$ associated to the choice of the Fourier transform $\Fa$. Sonin's space corresponds to the conditions $p^2-\lambda^2\geq 0$ and $q^2-\lambda^2\geq 0$ and the region of interest for the counting of eigenvalues is thus 
\begin{figure}
\begin{center}
\includegraphics[scale=0.8]{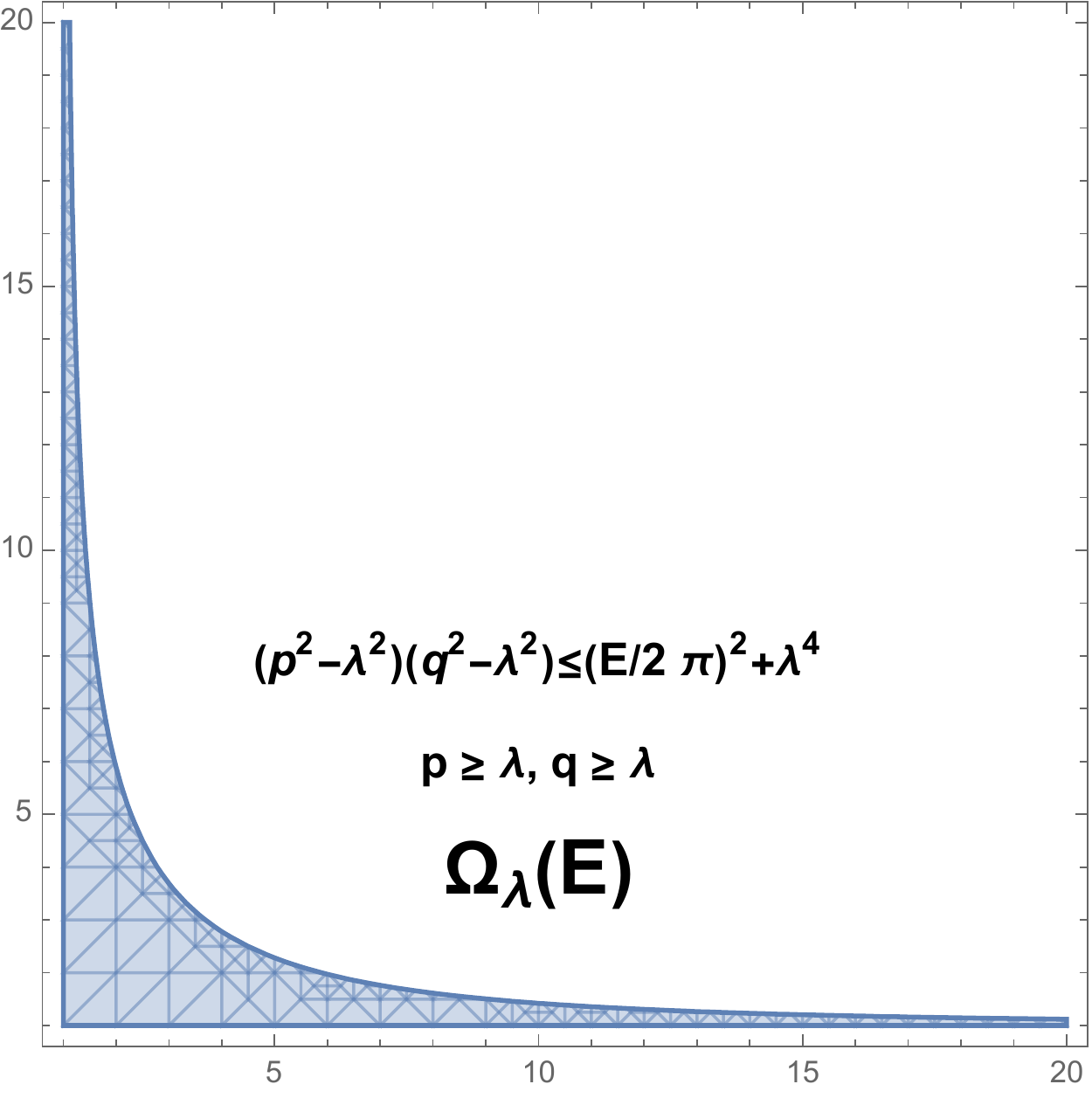}
\end{center}
\caption{The subset $\Omega_\lambda(E)$ in the first quadrant\label{figquant3}}
\end{figure}

$$
\Omega_\lambda(E):=\{(q,p)\mid q\geq \lambda, p\geq \lambda, H_\lambda(p,q)\leq \left(\frac{E}{2\pi}\right)^2+ \lambda^4\}$$ 
The area  of $\Omega_\lambda(E)$
 is given, with $a=\left(\frac{E}{2\pi}\right)^2+ \lambda^4$,  by the convergent integral
$$
I_\lambda(a)=\int_\lambda^{\infty}\left(\frac{\sqrt{a+\lambda^2 x^2-\lambda^4}}{\sqrt{x^2-\lambda^2}}-\lambda\right)dx
$$
One has, with $x=\lambda y$, the equality 
$$
I_\lambda(a)=\lambda\int_1^{\infty}\left(\frac{\sqrt{a+\lambda^4 y^2-\lambda^4}}{\sqrt{\lambda^2y^2-\lambda^2}}-\lambda\right)dy=\lambda^2
\int_1^{\infty}\left(\frac{\sqrt{a\,\lambda^{-4}+ y^2-1}}{\sqrt{y^2-1}}-1\right)dy
$$
Thus one obtains the equality
\begin{equation}\label{scalingrelation}
I_\lambda(a)=\lambda^2I_1(a\,\lambda^{-4})	
\end{equation}
We recall that the elliptic integrals $E(m)$ and $K(m)$ are defined by
$$
E(m):=\int_0^{\pi/2 } \sqrt{1-m \  \sin ^2 \theta} \, d\theta, \ \ 
K(m):=\int_0^{\pi/2 }\frac{1}{ \sqrt{1-m \  \sin ^2 \theta} }\, d\theta
$$

 \begin{lem}\label{integral} The integral $I(a)=I_1(a)$ is given by the sum of elliptic integrals
 \begin{equation}\label{integral1}
 	I(a)=a K(1-a)-E(1-a)+1
 \end{equation}	
\end{lem}
\proof One has, with $m=1-a$, $x=1/t$ 
$$
I(a)=\int_0^{1 }\left(\frac{\sqrt{1-m t^2 }}{\sqrt{1- t^2 }}-1\right)\frac{dt}{t^2}
$$
Let
$$
g(t):=-\frac{\sqrt{1-t^2} \left(m t^2+\sqrt{1-t^2} \sqrt{1-m t^2}-1\right)}{t \sqrt{1-m t^2}}
$$
One has $g(0)=0$, $g(1)=0$ and the derivative of $g$ is equal to 
$$
g'(t)=-\left(\frac{\sqrt{1-m t^2}}{\sqrt{1-t^2}}-1\right){t^{-2}}+\frac{1-m}{\sqrt{1-t^2} \sqrt{1-m t^2}}-\frac{\sqrt{1-m t^2}}{\sqrt{1-t^2}}+1
$$
so that the equality $\int_0^1g'(t)dt=0$ gives \eqref{integral1}.\endproof 
We thus get 
\begin{prop}\label{asympt1} The semiclassical approximation to the number of negative eigenvalues $\xi$  of $\wsa$ with $-\xi\leq E^2$ on even functions is the same as on odd functions and is equal to $2 \sigma(E,\lambda)$ where 
 \begin{equation}\label{sigmaE}
\sigma(E,\lambda)\sim \frac{E}{2 \pi }\left(\log \left(\frac{E}{2 \pi }\right)-1+\log (4)-2 \log (\lambda )\right)
+\lambda^2+o(1)
\end{equation}
 \end{prop}
\proof The semiclassical approximation corresponds, for the restriction to even functions (or to odd functions), to twice the area of $\Omega_\lambda(E)$ and hence to   
$
I_\lambda(a)=\lambda^2I(a\,\lambda^{-4})
$, for $a=\left(\frac{E}{2\pi}\right)^2+\lambda^4$. One has the asymptotic expansion for $a\to \infty$
\begin{equation}\label{asympt}
I(a)\sim \frac{1}{2} \sqrt{a} \left(\log (a)-2+4 \log (2)\right)+1+\frac{1}{8} \sqrt{\frac{1}{a}} (-\log (a)-4 \log (2))+0 \left(\frac{1}{a}\right) 
\end{equation}
so that
\begin{equation}\label{asympt}
I_\lambda(a)\sim \frac{1}{2} \sqrt{a} \left(\log (a)-2+4 \log (2)-4\log \lambda\right)+\lambda^2+o(1)
\end{equation}
  We then use the expansions 
$$
\sqrt{a}=\frac{E}{2\pi}+ O(1/E), \ \ \log(a)=2 \log \left(\frac{E}{2\pi}\right)+O(1/E^2)
$$
and obtain \eqref{sigmaE}.\endproof 

\section{Dirac operators}\label{sectdirac}
The results of Section \ref{sectsemiclassical} show that for suitable values of $\lambda$ the negative spectrum of $\wsa$ has the same ultraviolet behavior as the squares of zeros of the Riemann zeta function. Since $\wsa$ is a differential operator of second order we liken it to the Klein-Gordon operator and construct the analogue of the Dirac operator. We first use the Darboux process (see \cite{Deift}, \cite{grunbaum2}) to factorize $\wsa$ as a product of two first order differential operators.
\begin{lem}\label{darboux1} Let $p(x)=x^2-\lambda^2$, $V(x)=4\pi^2 \lambda^2 x^2$,  $L=\partial(p(x)\partial)+V(x)$, $(\nabla f)(x):=p(x)^{1/2}\partial f(x)$ and $U$ the unitary operator 
$$
 U: L^2([\lambda,\infty),dx)\to L^2([\lambda,\infty),p(x)^{-1/2}dx), \ \ U(\xi)(x):=p(x)^{1/4}\xi(x).
$$
Let $w(x)$ be a solution of the  equation 
\begin{equation}\label{Ricca0}
\nabla w(x)+w(x)^2=-V(x)+\left(\frac{p''(x)}{4}-\frac{p'(x)^2}{16 p(x)}\right)\qqq x\in [\lambda,\infty)
\end{equation}  then one has $L=U^*(\nabla+w)(\nabla-w)U$.
  \end{lem} 
  \proof Let $f$ be a smooth function on $\R$ and consider the differential operators $T_1:=f\partial_x f $ and  $T_2:=\partial_x f^4\partial_x  $. Let us show that $T_1^2-T_2$ is an operator of order zero: one has
  $$
  T_1^2=f\partial_x f^2\partial_x f=-f'f^2\partial_x f+\partial_x f^3\partial_x f
  $$ 
  $$
 -f'f^2\partial_x f=-f'^2f^2-f'f^3\partial_x, \ \ \partial_x f^3\partial_x f=\partial_x f^4\partial_x+\partial_x f^3 f'
  $$
  so that $T_1^2-T_2$ is the multiplication by $2f'^2f^2+f^3f''$. Applying this for $f(x)=p(x)^{1/4}$ gives 
  $$
  (U^*\nabla U)^2=\partial_x p(x)\partial_x+\frac{p''(x)}{4}-\frac{p'(x)^2}{16 p(x)}
  $$
  from which the conclusion follows using \eqref{Ricca0}. \endproof

  We now determine all solutions of the Riccati equation \eqref{Ricca0} which gives
\begin{equation}\label{Ricca1}
\sqrt{x^2-\lambda^2}\, w'(x)+w(x)^2=-4\pi^2\lambda^2 x^2-\frac{1}{4}\ \frac{x^2}{  x^2-\lambda^2}+\frac{1}{2}\end{equation}
The next Lemma is standard using the reduction of a Riccati equation to a Bernoulli equation.
\begin{lem}\label{liouville2dirac2} Let $u_j$ be  two real valued solutions of $Lu=0$ which generate  the linear space of solutions in $(\lambda,\infty)$. \newline
$(i)$~For $z\in \C$ and $u=u_1+z u_2$ the solution $u$ has no zero in $(\lambda,\infty)$ if $z \notin \R$ and an infinity of zeros otherwise. \newline
$(ii)$~All solutions of the Riccati equation \eqref{Ricca1} are given by 
\begin{equation}\label{anzats1}
w_z(x)=\frac{(x^2-\lambda^2)^{1/4} \partial \left((x^2-\lambda^2)^{1/4} u(x)\right)}{u(x)}
\end{equation} where $u=u_1+z u_2$ and $z\in \C\setminus \R$.\newline
$(iii)$~The map $z\mapsto w_z$ from $\C\setminus \R$ to the space of solutions of \eqref{Ricca1} is a homeomorphism.	
\end{lem}
\proof $(i)$~Let $x\in (\lambda,\infty)$ and $z\in \C\setminus \R$, $u=u_1+z u_2$. Assume $u(x)=0$. Then $u_1(x)+z u_2(x)=0$ and since $z\in \C\setminus \R$ this implies $u_1(x)=0$ and $u_2(x)=0$. The Wronskian $p(x)(u_1'(x)u_2(x)-u_2'(x)u_1(x))$ is constant and non-zero since the $u_j$ are independent solutions. Thus we get a contradiction and $u$ has no zero. Moreover since the equation $Lu=0$ is in the LCO case any real valued solution has infinitely many zeros.  \newline
$(ii)$~The standard solution of the Riccati equation \eqref{Ricca0} is of the form
\begin{equation}\label{anzats}
w(x)=\frac{p(x)^{1/4} \partial \left(p(x)^{1/4} u(x)\right)}{u(x)}
\end{equation}
which gives 
$$
\nabla w(x)+w(x)^2=\frac{p''(x)}{4}+\frac{p'(x) u'(x)}{u(x)}-\frac{p'(x)^2}{16 p(x)}+\frac{p(x) u''(x)}{u(x)}
$$
so that 
$$
Lu=0\Rightarrow \nabla w(x)+w(x)^2=-V(x)+\left(\frac{p''(x)}{4}-\frac{p'(x)^2}{16 p(x)}\right)
$$
Thus by $(i)$  any $w_z$, $z\in \C\setminus \R$ is a solution of 
\eqref{Ricca1}. Using the three values $\{i,-i,j\}$ for $z$ and the reduction to a Bernoulli equation one can express the general solution of 
\eqref{Ricca1} in the form 
\begin{equation}\label{anzats2.5}
w=w_i+\frac{(w_{-i}-w_i)(w_{j}-w_i)}{(1-t)(w_{j}-w_i)+t\,(w_{-i}-w_i) }=w_{z(t)}
\end{equation}
where 
$$
z(t)=\frac{i (i (t-1)+j (t+1))}{i (t+1)+j (t-1)}\in \C\setminus \R
$$

$(iii)$~The formula \eqref{anzats2.5} establishes an homeomorphism between the space of solutions of the Riccati equation and the complement in $\P^1(\C)$ of the circle 
$$
\{t\in \P^1(\C)\mid z(t)\in \P^1(\R)\}
$$
and thus the map $z\mapsto w_z$ from $\C\setminus \R$ to the space of solutions of \eqref{Ricca1} is a homeomorphism.\endproof 

\begin{prop}\label{dirac1} Let $w$ be a  solution of the Riccati equation \eqref{Ricca1}
and $\Dirac$ be the matrix of order one operators 
 	\begin{equation}\label{sqrt}
 \Dirac=	\left(
\begin{array}{cc}
 0 & \nabla+w(x) \\
\nabla-w(x) & 0 \\
\end{array}
\right)
 	\end{equation}
Then the square of $\Dirac$ is diagonal with each diagonal term spectrally equivalent to $L$,
$$
U^*\Dirac^2 U=\left(
\begin{array}{cc}
 L & 0 \\
0 & L+2\nabla w(x) \\
\end{array}
\right)
$$
\end{prop}
 The proof is straightforward. The use of the Darboux process in this construction is related to the theory of isospectral deformations  \cite{Deift, grunbaum2}.
 \section{Ultraviolet behavior of spectrum of Dirac, case $\lambda=\sqrt 2$}\label{sectroot2}
 
In this section we take  $\lambda=\sqrt 2$, and consider the operator 
$2 \Dirac$ where $\Dirac$ is as defined in  Proposition \ref{dirac1}. 
\begin{thm}\label{thm2D} The operator $2 \Dirac$ has discrete simple spectrum contained in $\R\cup i\R$. Its  imaginary eigenvalues are symmetric under complex conjugation and the counting function $N(E)$ counting those of positive imaginary part less than $E$ fulfills
\begin{equation}\label{ne}
N(E)\sim \frac{E}{2 \pi }\left(\log \left(\frac{E}{2 \pi }\right)-1\right)+O(1)
\end{equation}
\end{thm}
\proof By  Proposition \ref{dirac1} the spectrum of $2\Dirac$ consists of the complex numbers of the form $\xi=\pm 2 \sqrt \alpha$ where $\alpha$ varies in the spectrum of $L$. The latter is real and the number of negative eigenvalues $\alpha\geq -E^2$ is given by Proposition \ref{asympt1} as $2 \sigma(E,\lambda)$, thus selecting the root with positive imaginary part one gets 
$$
0<\Im(\xi)\leq E \iff \alpha\geq -(E/2)^2
$$ 
and the number $N(E)$ of such $\xi$ is thus
$$
2 \sigma(E/2,\sqrt 2)= \frac{E}{2 \pi }\left(\log \left(\frac{E/2}{2 \pi }\right)-1+\log (4)-2 \log (\sqrt 2 )\right)+O(1)=
$$
$$
=\frac{E}{2 \pi }\left(\log \left(\frac{E}{2 \pi }\right)-1\right)+O(1)
$$
 which gives the required estimate.\endproof 

 \begin{figure}[H]	\begin{center}
\includegraphics[scale=0.8]{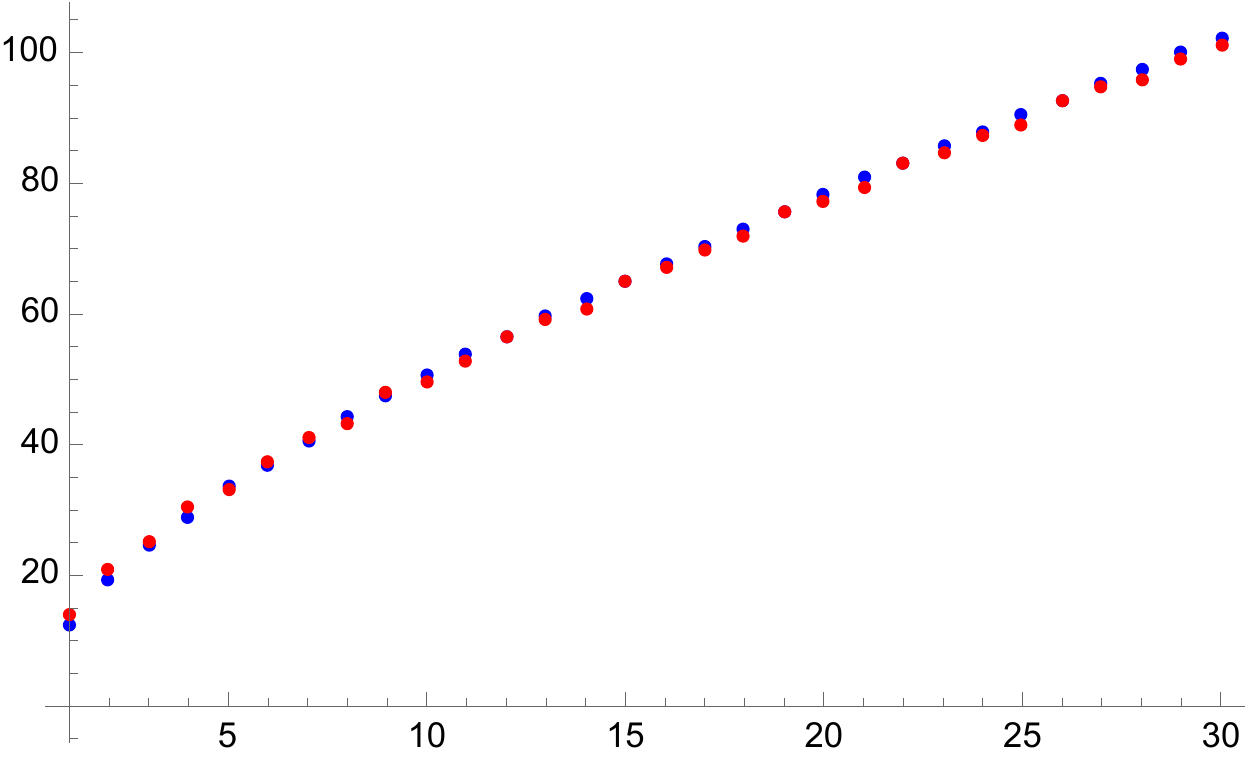}
\end{center}
\caption{The spectrum (in blue) compared to zeros of zeta (red)\label{figquant3}}
\end{figure}
 \begin{figure}[H]	\begin{center}
\includegraphics[scale=0.8]{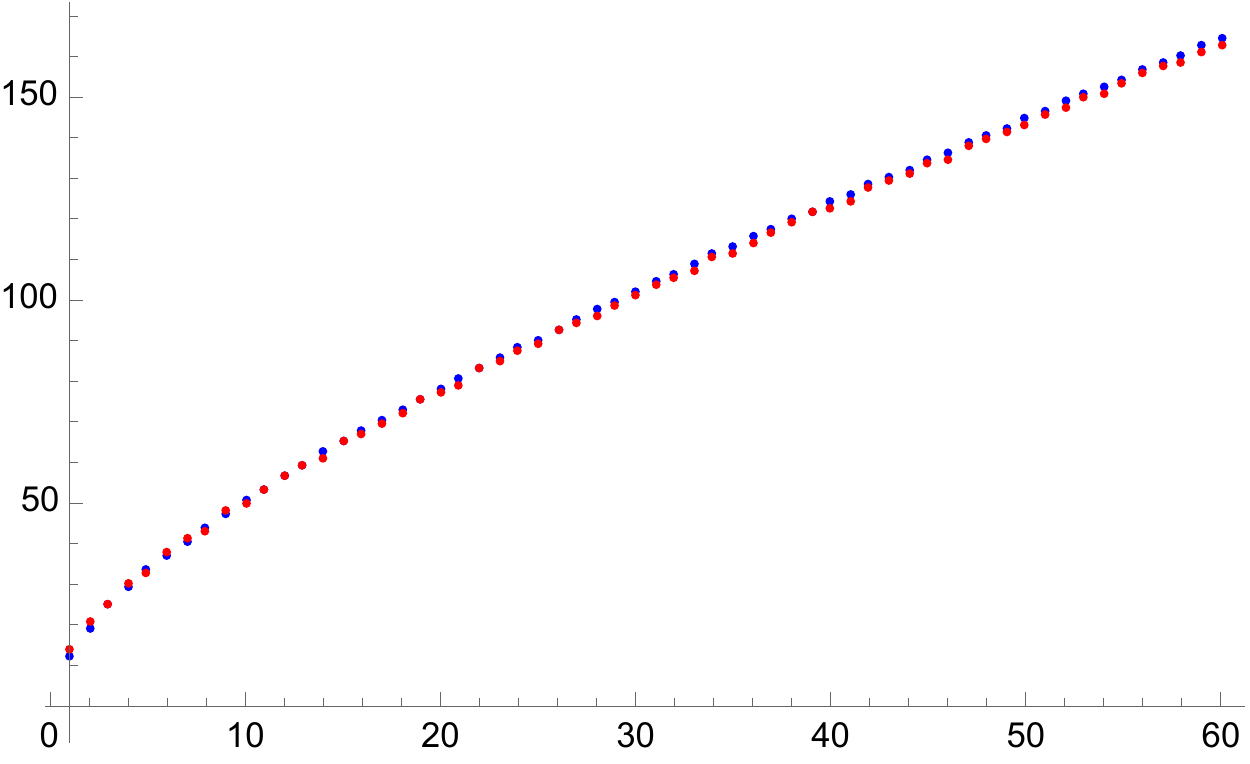}
\end{center}
\caption{The spectrum compared to zeros of zeta\label{figquant3}}
\end{figure} \begin{figure}[H]	\begin{center}
\includegraphics[scale=0.8]{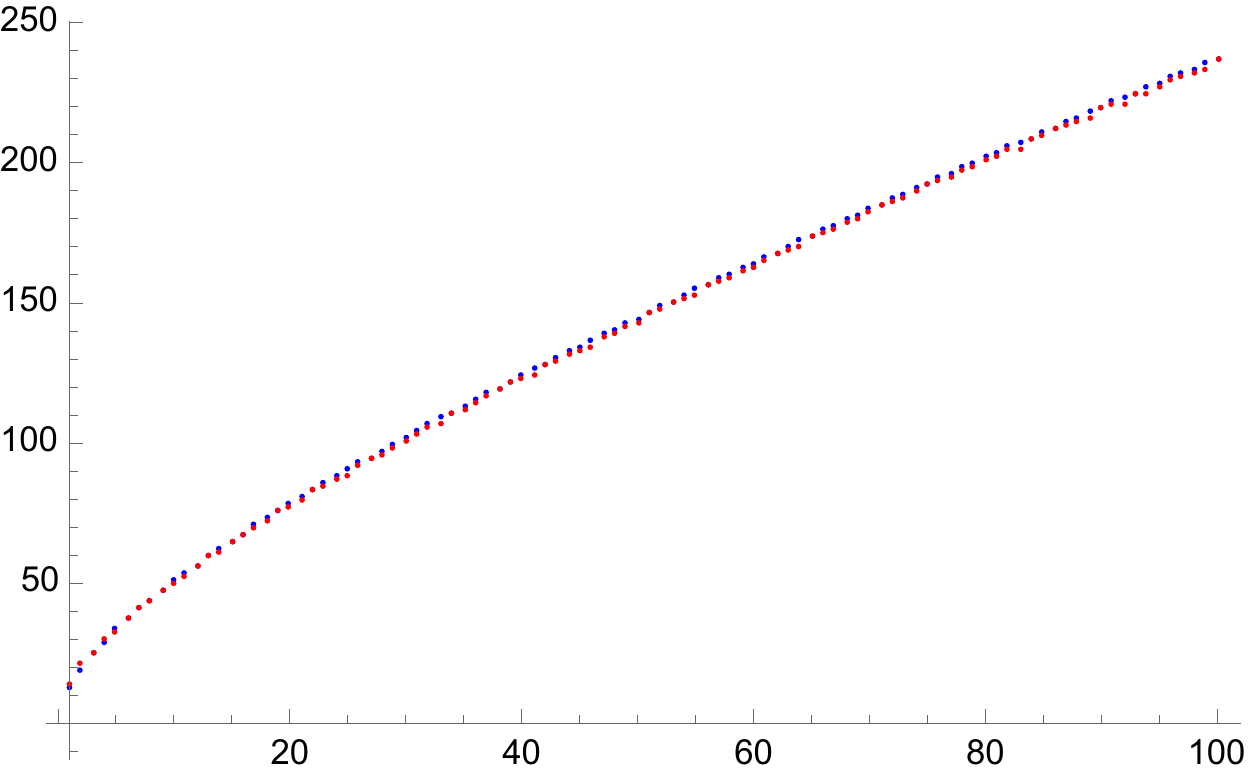}
\end{center}
\caption{The spectrum compared to zeros of zeta\label{figquant3}}
\end{figure}
 
\section{Final remarks}\label{sectrems}

We gather in this final section a number of more speculative remarks.

\subsection{Geometric meaning of Theorem \ref{thm2D}}

The operator $2\Dirac$ of Theorem \ref{thm2D} together with the action by multiplication of smooth functions on the interval $[\sqrt 2,\infty)$ would define a spectral triple if 
$2\Dirac$ were self-adjoint (or skew adjoint) but its spectrum contains both real and imaginary pieces. The leading term $\nabla=\left(2\sqrt{x^2-2}\right)\partial_x$ shows that the 
corresponding classical metric is (with $\lambda=\sqrt 2$ from now on)
$$
ds^2=-\frac 14 dx^2/(x^2-2)=\frac{1}{\alpha(x)}dx^2, \ \ \alpha(x)=-4 (x^2-2)
$$
This $ds^2$ changes sign when crossing the boundary $x=\sqrt 2$ and this suggests,
in order to handle all even functions on $\R$ and to take into account the real and imaginary eigenvalues of the square of $2\Dirac$, to 
 look for a two dimensional metric with signature $(-1,1)$ of the form
$$
ds^2=-\alpha(x)dt^2+\frac{1}{\alpha(x)}dx^2
$$
This geometry corresponds to a black hole in two space-time dimensions with horizon at $x=\pm \sqrt 2$. It fulfils the 2-dimensional analogue of Einstein's equation with a cosmological constant $=8$ and no source \cite{MST}. One can look at the null curves and this means 
$$
\frac{dx}{dt}=\alpha(x)=-4 (x^2-2) \Rightarrow t(x)=\frac {1}{8\sqrt 2 }  \log \left((\sqrt 2 +x)/(x-\sqrt 2 )\right)+c
$$
one then passes to the new coordinates $v=t-t(x)$, and $x$ unchanged. In these new coordinates one re-expresses the metric in the smooth form 
$$
ds^2=4 \left(x^2-2\right)dv^2-2 dv dx
$$ 
\begin{figure}[H]	\begin{center}
\includegraphics[scale=1.02]{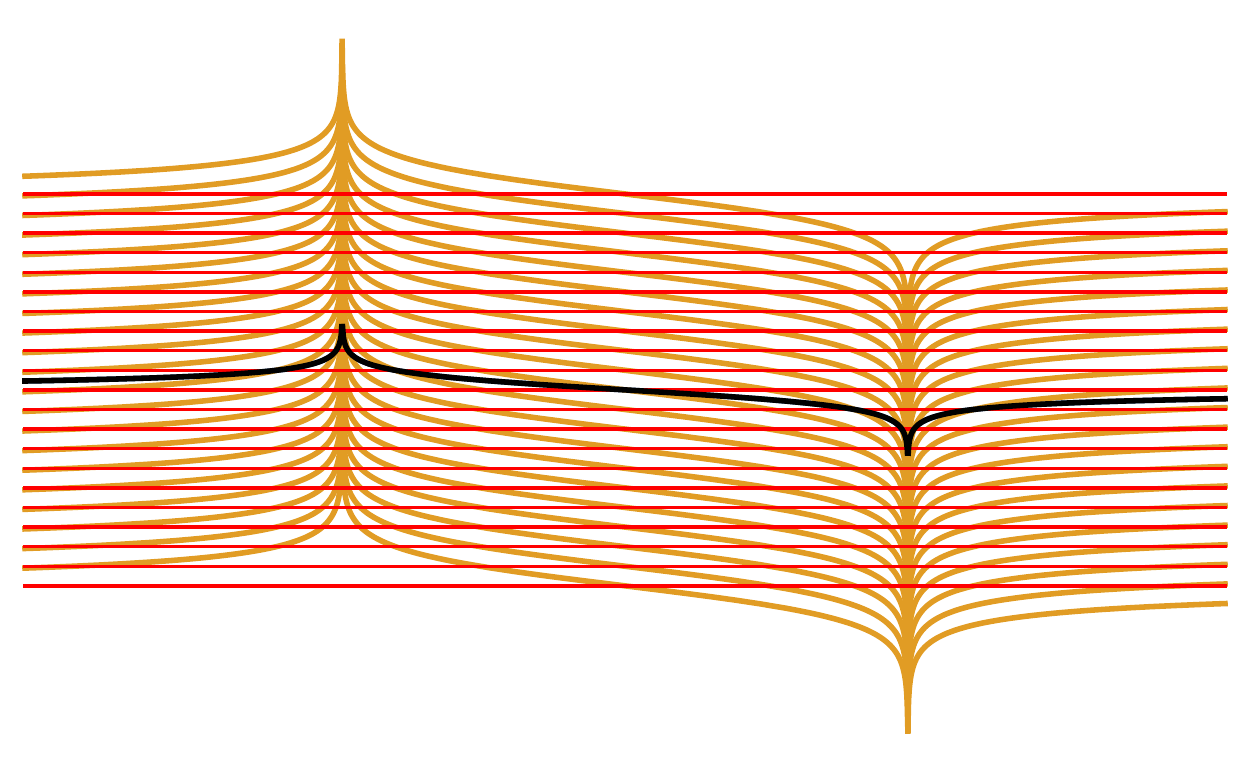}
\end{center}
\caption{Light rays in the two dimensional geometry, and in black the original curve. The vertical lines are the horizons at $x=\pm \sqrt 2$.\label{licorne}}
\end{figure}

In this metric the light rays are given by $v=v_0$ (\ie the horizontal lines of Figure \ref{licorne}) and by the curves 
$$
v(x)=\frac{1}{4 \sqrt{2}}\log\bigg\vert \frac{x-\sqrt{2}}{x+\sqrt{2}}\bigg\vert+c
$$
\ie the solutions of the equation 
$
dv=\frac{dx}{2\left(x^2-2\right)}
$.
The original curve given by $t=0$ corresponds to the graph of $v=-t(x)$ as shown in black in Figure \ref{licorne}.
\subsection{Positive eigenvalues of $\wsa$ and trivial zeros of Zeta} \label{sectposeigen}
 The eigenvalues $\chi(n)$ of the restriction of $\wsa$ to even functions in the interval $[-\lambda,\lambda]$ have a well understood asymptotic form which by \cite{ORX} Theorem 3.11 implies that, independently of the value of $\lambda$, (note that we only consider even functions so that the index $n$ of \opcit is replaced by $2n$) 
 $$
 \chi(n)= \left(2n+\frac 12\right)^2+ O(1), \ n\to \infty
 $$  
  This behavior is the same as that of the squares of the trivial zeros of the Riemann zeta function with the same shift of $\frac 12$ as for the critical line. To obtain a convincing relation one would need to analyze the extension of $2\Dirac$ to (two copies of) the even functions on $\R$ as well as the conditioning  of the Hilbert space needed to eliminate the positive square roots of the $ \chi(n)$. 

\subsection{Spectral truncation}

In order to eliminate the real eigenvalues of $2 \Dirac$ coming from the positive eigenvalues of $\wsa$ one can effect a spectral truncation \cite{CS1},  the algebra of functions acting by multiplication is then replaced by the operator system obtained by compression on Sonin's space. In a similar manner one can use spectral truncation to eliminate the positive square roots with the notations of \S \ref{sectposeigen}.

\end{document}